\documentclass[twoside,12pt,reqno]{amsart}
\usepackage{latexsym}%%triangles:\lhd\rhd\unlhd\unrhd%
\newcommand{\qdn}{\hspace*{-1.5mm}}
\newcommand{\qqdn}{\hspace*{-2.5mm}}
\newcommand{\xqdn}{\hspace*{-5.0mm}}
\newcommand{\xxqdn}{\hspace*{-10mm}}

\newcommand{\fns}{\footnotesize}

\newcommand{\ffnk}[4]{\left[\qdn\ba{#1}#3\\#4\ea{\!\Big|\:#2}\right]}

\newcommand{\binm}{\binom}

\newcommand{\nnm}{\nonumber}
\newcommand{\be}{\begin{equation}}
\newcommand{\ee}{\end{equation}}
\newcommand{\ba}{\begin{array}}
\newcommand{\ea}{\end{array}}
\newcommand{\bmn}{\begin{eqnarray}}
\newcommand{\emn}{\end{eqnarray}}
\newcommand{\bnm}{\begin{eqnarray*}}
\newcommand{\enm}{\end{eqnarray*}}
\newcommand{\bln}{\begin{subequations}}
\newcommand{\eln}{\end{subequations}}

\newtheorem{thm}{Theorem}%[section]
\newtheorem{lemm}[thm]{Lemma}

\newtheorem{entry}{Entry}%%%%%%%%%%%%%%%%

\newcommand{\bbtm}[4]{\bibitem{kn:#1}{#2,}~{#3,}~{#4.}}
\newcommand{\cito}[1]{\cite{kn:#1}}
\newcommand{\citu}[2]{\cite[#2]{kn:#1}}

\begin{document} %%%%%%%%%% This paper is published in %%%%%%%
{\fns% \today\hfill\copyright%% Printed in China} %%%%%%%%%%%%%%%
%%%%%%%%%%%%%%%%%%%%%%%%%%%%%%%%%%%%%%%%%%%%%%%%%%%%%%%%%%%%%%
\title{The united proofs for three $q$-extensions of Dougall's $_2H_2$ summation formula}
\author{Chuanan Wei}


\dedicatory{
Department of Medical Informatics\\
  Hainan Medical University, Haikou 571199, China}
\thanks{\emph{Email address}:
      weichuanan78@163.com}

\address{ }

\begin{abstract}
In terms of the analytic continuation method, we give the united
proofs for three $q$-extensions of Dougall's $_2H_2$ summation
formula. Some related results are also discussed in this paper.
\end{abstract}

%%%%%%%%%%%%%%%%%%%%%%%%%%%%%%%%%%%%%%%%%%%%%%%%%%%%%%%%%%%%%%%%%%%
\maketitle\thispagestyle{empty}%%%%%%%%%%%%%%%%%%%%%%%%%%%%%%%%%%%%
\markboth{Chuanan Wei}%%%%%%%%%%%%%%%%%%%%%%%%%%%%
         {$q$-Extensions of Dougall's $_2H_2$ summation formula}

%%%%%%%%%%%%%%%%%%%%%%%%%%%%%%%%%%%%%%%%%%%%%%%%%%%%%%%%%%%%%%%%%%%
%%%%%%%%%%%%%%%%%%%%%%%%%%%%%%%%%%%%%%%%%%%%%%%%%%%%%%%%%%%%%%%%%%%
\section{Introduction}
%%%%%%%%%%%%%%%%%%%%%%%%%%%%%%%%%%%%%%%%%%%%%%%%%%%%%%%%%%%%%%%%%%%
For a complex number $x$, define the gamma function by Euler's
integral
\[\Gamma(x)=\int_{0}^{\infty}t^{x-1}e^{-t}dt\quad\text{with}\quad\mathfrak{Re}(x)>0.\]
Then the shifted factorial can be expressed as
\[(x)_n=\frac{\Gamma(x+n)}{\Gamma(x)},\]
 where $n$ is an arbitrary integer.

Following Slater~\cito{slater}, define the bilateral hypergeometric
series to be
 \bnm
&&_{r}H_s\ffnk{cccc}{z}{a_1,&a_2,&\cdots,&a_r}{b_1,&b_2,&\cdots,&b_s}
 \:=\:\sum_{k=-\infty}^\infty
\frac{(a_1)_k(a_2)_k\cdots(a_r)_k}{(b_1)_k(b_2)_k\cdots(b_s)_k}z^k.
 \enm
In 1907, Dougall \cito{dougall} derived the beautiful identity
\bmn\label{dougall}
 {_2H_2}\ffnk{ccc}{1}{a,\:b}{c,\:d}
=\frac{\Gamma(1-a)\Gamma(1-b)\Gamma(c)\Gamma(d)\Gamma(c+d-a-b-1)}{\Gamma(c-a)\Gamma(c-b)\Gamma(d-a)\Gamma(d-b)},
 \emn
 provided $\mathfrak{Re}(c+d-a-b)>1$, according to the contour integral method. Different proofs of it can  be
 found in  \citu{andrews-a}{Section 2.8}, \cito{chu-a}, \cito{schlosser},
\citu{slater}{Section 6.1}, and \cito{wei}.

Subsequently, define the $q$-shifted factorial by
 \bnm
(x;q)_{\infty}=\prod_{i=0}^{\infty}(1-xq^i),\quad
(x;q)_n=\frac{(x;q)_{\infty}}{(xq^n;q)_{\infty}},
 \enm
where $x$, $q$  are complex numbers satisfying the condition $|q|<1$
and $n$ is an arbitrary integer. For convenience, we shall adopt the
two simplified notations:
 \bnm
&&(x_1,x_2,\cdots,x_r;q)_{\infty}=(x_1;q)_{\infty}(x_2;q)_{\infty}\cdots(x_r;q)_{\infty},\\
&&(x_1,x_2,\cdots,x_r;q)_{n}=(x_1;q)_{n}(x_2;q)_{n}\cdots(x_r;q)_{n}.
 \enm
Following Gasper and Rahman \cito{gasper}, define the bilateral
basic hypergeometric series to be
 \bnm
\quad{_r\psi_s}\ffnk{cccccc}{q;z}{a_1,&a_2,&\cdots,&a_r}{b_1,&b_2,&\cdots,&b_s}
  =\sum_{k=-\infty}^{\infty}\frac{(a_1,a_2,\cdots,a_r;q)_k}{(b_1,b_2,\cdots,b_s;q)_k}
\Big\{(-1)^kq^{\binm{k}{2}}\Big\}^{s-r}z^k.
 \enm
Thus Ramanujan's $_1\psi_1$ summation formula (cf.
\citu{gasper}{Appendix II.29}) can be stated as
 \bmn\label{ramanujan}
{_1\psi_1}\ffnk{cccccc}{q;z}{a}{c}=
\frac{(q,c/a,az,q/az;q)_{\infty}}{(c,q/a,z,c/az;q)_{\infty}},
 \emn
provided $|c/a|<|z|<1$. Equation \eqref{ramanujan} is very important
in the theory of special functions. Several beautiful proofs of it
can be enjoyed in the papers
\cite{kn:andrews-b,kn:andrews-c,kn:andrews-d,kn:chen-b,kn:ismail}.
The case $b_1=q$ of the bilateral basic hypergeometric series is
exactly the unilateral basic hypergeometric series
 \bnm\quad
{_{r}\phi_{s-1}}\ffnk{cccccc}{q;z}{a_1,&a_2,&\cdots,&a_r}{&b_2,&\cdots,&b_s}
  =\sum_{k=0}^{\infty}\frac{(a_1,a_2,\cdots,a_r;q)_k}{(q,b_2,\cdots,b_s;q)_k}
\Big\{(-1)^kq^{\binm{k}{2}}\Big\}^{s-r}z^k.
 \enm
Then three different $q$-extensions of \eqref{dougall} can be laid
out as follows.

\begin{thm} \label{thm-a}
Let $a$, $b$, $c$, $d$ be complex numbers. Then
 \bnm
\qdn{_2\psi_2}\ffnk{cccccc}{q;z}{a,b}{c,d}=\frac{(az,c/b,d/a,qd/abz;q)_{\infty}}{(z,d,q/b,cd/abz;q)_{\infty}}
{_2\psi_2}\ffnk{cccccc}{q;\frac{d}{a}}{a,abz/d}{c,az},
 \enm
where $\max\{|z|, |cd/abz|, |d/a|, |c/b|\}<1$.
\end{thm}

\begin{thm} \label{thm-b}
Let $a$, $b$, $c$, $d$ be complex numbers. Then
 \bnm
{_2\psi_2}\ffnk{cccccc}{q;z}{a,b}{c,d}
&&\xqdn\!=\frac{(q,b,c/a,d/a,az,q/az;q)_{\infty}}{(c,d,q/a,b/a,z,q/z;q)_{\infty}}
 {_2\phi_1}\ffnk{cccccc}{q;\frac{cd}{abz}}{qa/c,qa/d}{qa/b}\\
&&\xqdn\!+\:\:idem(a;b),
 \enm
where the convergent condition is $\max\{|z|, |cd/abz|\}<1$ and the
symbol $idem(a;b)$ after an expression signifies that the front
expression is repeated with a and b interchanged.
\end{thm}

\begin{thm} \label{thm-c}
Let $a$, $b$, $c$, $d$ be complex numbers. Then
 \bnm
{_2\psi_2}\ffnk{cccccc}{q;z}{a,b}{c,d} &&\xqdn=\:
\frac{(q,c/b,q/d,abz/d,qd/abz;q)_{\infty}}{(q/a,q/b,c,az/d,cd/abz;q)_{\infty}}
{_2\phi_1}\ffnk{cccccc}{q;\frac{qb}{d}}{cd/abz,d/a}{qd/az}
 \nnm\\&&\xqdn
 -\,\,\frac{(q,b,q/d,qc/d,d/a,az/q,q^2/az;q)_{\infty}}{(q/a,c,d/q,q^2/d,qb/d,az/d,qd/az;q)_{\infty}}
 {_2\phi_1}\ffnk{cccccc}{q;z}{qa/d,qb/d}{qc/d},
 \enm
provided $\max\big\{|z|,|cd/abz|,|qb/d|\big\}<1$.
\end{thm}

In 1950, Bailey \cito{bailey} established Theorem 1 by applying the
method of comparing coefficients to the product of
\eqref{ramanujan}.
 For the semi-finite form of it, the reader
is referred to Chen and Fu \cito{chen-b}. Gasper and Rahman
\citu{gasper}{Section 5.4} have shown that there exist two
expansions of an $_r\psi_r$ series by means of $r$ $_r\phi_{r-1}$
series (cf. \citu{gasper}{Equations (5.4.4) and (5.4.5)}). The case
$r=2$ of the former is exactly Theorem \ref{thm-b}. when $d=b$, it
becomes
 \bnm
{_1\psi_1}\ffnk{cccccc}{q;z}{a}{c}=
\frac{(q,c/a,az,q/az;q)_{\infty}}{(c,q/a,z,q/z;q)_{\infty}}
 {_1\phi_0}\ffnk{cccccc}{q;\frac{c}{az}}{qa/c}{-}.
 \enm
Evaluating the $_1\phi_0$ series on the right hand side by
$q$-binomial theorem (cf. \citu{gasper}{Appendix II.3}):
 \bnm
 {_1\phi_0}\ffnk{cccccc}{q;z}{a}{-}=\frac{(az;q)_{\infty}}{(z;q)_{\infty}},
 \enm
 we get Ramanujan's $_1\psi_1$ summation formula \eqref{ramanujan}. Theorem \ref{thm-c} was deduced by Chens and Gu \cito{chen-a}
  in accordance with Cauchy's method. Interestingly,
this theorem reduces directly to \eqref{ramanujan} when $d=a$.
Recently, the research of $q$-congruence associated with summation
and transformation formulas for the unilateral basic hypergeometric
series attracts several mathematicians. Some nice results can be
seen in the papers \cite{kn:guo-a,kn:guo-b,kn:guo-c,kn:guo-d}.

 A property of the analytic function (cf. \citu{lang}{p.90}; see also \cite{kn:askey,kn:ismail,kn:zhu}), which plays a central role in this paper,
 can be displayed as the following lemma.

\begin{lemm} \label{lemma}
Let $U$ be a connected open set and $f$, $g$ be analytic on $U$. If
$f$ and $g$ agree infinitely often near an interior point of $U$,
then we have $f(z)=g(z)$ for all $z\in U$.
\end{lemm}

The structure of the paper is arranged as follows. By the
utilization of Lemma \ref{lemma}, we shall supplies the united
proofs
 of Theorems \ref{thm-a}-\ref{thm-c} in Section 2.
Some related results are also discussed in Section 3.
%%%%%%%%%%%%%%%%%%%%%%%%%%%%%%%%%%%%%%%%%%%%%%%%%%%%%%%%%%%%%%%%%%%
\section{Proofs of Theorems \ref{thm-a}-\ref{thm-c}}
%%%%%%%%%%%%%%%%%%%%%%%%%%%%%%%%%%%%%%%%%%%%%%%%%%%%%%%%%%%%%%%%%%%
%%%%%%%%%%%%%%%%%%%%%%%%%%%%%%%%%%%%%%%%%%%%%%%%%%%%%%%%%%%%%%%%%%%
\subsection{Proof of Theorem \ref{thm-a}}
%%%%%%%%%%%%%%%%%%%%%%%%%%%%%%%%%%%%%%%%%%%%%%%%%%%%%%%%%%%%%%%%%%%
For a positive integer $m$, we have the relation
 \bmn\label{bailey-a}
{_2\psi_2}\ffnk{cccccc}{q;z}{a,b}{q^{1+m},d}
&&\xqdn=\sum_{k=-m}^{\infty}\frac{(a,b;q)_k}{(q^{1+m},d;q)_k}z^k
=\sum_{k=0}^{\infty}\frac{(a,b;q)_{k-m}}{(q^{1+m},d;q)_{k-m}}z^{k-m}
 \nnm\\
&&\xqdn=\frac{(a,b;q)_{-m}z^{-m}}{(q^{1+m},d;q)_{-m}}{_2\phi_1}\ffnk{cccccc}{q;z}{aq^{-m},bq^{-m}}{dq^{-m}}.
 \emn
Using \eqref{bailey-a} and Heine's transformation formula between
two $_2\phi_1$ series (cf. \citu{gasper}{Appendix III.2}):
\[{_2\phi_1}\ffnk{cccccc}{q;z}{a,b}{c}=\frac{(c/b,bz;q)_{\infty}}{(c,z;q)_{\infty}}{_2\phi_1}\ffnk{cccccc}{q;\frac{c}{b}}{abz/c,b}{bz},\]
we gain
  \bnm
&&\xqdn{_2\psi_2}\ffnk{cccccc}{q;z}{a,b}{q^{1+m},d}\\
&&\qqdn\:=\:\frac{(a,b;q)_{-m}z^{-m}}{(q^{1+m},d;q)_{-m}}\frac{(d/a,azq^{-m};q)_{\infty}}{(dq^{-m},z;q)_{\infty}}
{_2\phi_1}\ffnk{cccccc}{q;\frac{d}{a}}{abzq^{-m}/d,aq^{-m}}{azq^{-m}}\\
&&\qqdn\:=\:\frac{(a,b;q)_{-m}z^{-m}}{(q^{1+m},d;q)_{-m}}\frac{(d/a,azq^{-m};q)_{\infty}}{(dq^{-m},z;q)_{\infty}}
\qdn\sum_{k=-m}^{\infty}\qdn\frac{(abzq^{-m}/d,aq^{-m};q)_{k+m}}{(q,azq^{-m};q)_{k+m}}\bigg(\frac{d}{a}\bigg)^{k+m}\\
&&\qqdn\:=\:\frac{(az,q^{1+m}/b,d/a,qd/abz;q)_{\infty}}{(z,d,q/b,q^{1+m}d/abz;q)_{\infty}}
{_2\psi_2}\ffnk{cccccc}{q;\frac{d}{a}}{a,abz/d}{q^{1+m},az}.
 \enm
Split the ${_2\psi_2}$ series on both sides into two parts to
achieve
 \bmn\label{bailey-b}
&&\xqdn\sum_{k=0}^{\infty}\frac{(a,b;q)_k}{(q^{1+m},d;q)_k}z^k
+\sum_{k=1}^{\infty}\frac{(q/d;q)_k\prod_{i=1}^k(q^{1+m}-q^i)}{(q/a,q/b;q)_k}\bigg(\frac{d}{abz}\bigg)^k
 \nnm\\\nnm
&&\xqdn\:=\:\frac{(az,q^{1+m}/b,d/a,qd/abz;q)_{\infty}}{(z,d,q/b,q^{1+m}d/abz;q)_{\infty}}\\
&&\xqdn\:\times\:\bigg\{\sum_{k=0}^{\infty}\frac{(a,abz/d;q)_k}{(q^{1+m},az;q)_k}\bigg(\frac{d}{a}\bigg)^k
+\sum_{k=1}^{\infty}\frac{(q/az;q)_k\prod_{i=1}^k(q^{1+m}-q^i)}{(q/a,qd/abz;q)_k}\bigg(\frac{1}{b}\bigg)^k\bigg\}.
 \emn
Define two functions $f(c)$ and $g(c)$ by
 \bnm
&&f(c)=\sum_{k=0}^{\infty}\frac{(a,b;q)_k}{(c,d;q)_k}z^k
+\sum_{k=1}^{\infty}\frac{(q/d;q)_k\prod_{i=1}^k(c-q^i)}{(q/a,q/b;q)_k}\bigg(\frac{d}{abz}\bigg)^k,\\
&&g(c)=\frac{(az,c/b,d/a,qd/abz;q)_{\infty}}{(z,d,q/b,cd/abz;q)_{\infty}}\\
&&\qquad\times\:\bigg\{\sum_{k=0}^{\infty}\frac{(a,abz/d;q)_k}{(c,az;q)_k}\bigg(\frac{d}{a}\bigg)^k
+\sum_{k=1}^{\infty}\frac{(q/az;q)_k\prod_{i=1}^k(c-q^i)}{(q/a,qd/abz;q)_k}\bigg(\frac{1}{b}\bigg)^k\bigg\}.
 \enm
Then \eqref{bailey-b} shows that
 \bmn \label{bailey-c}
 f(c)=g(c)
 \emn
for $c=q^{1+m}$. According to Lemma \ref{lemma}, \eqref{bailey-c} is
correct for all $|c|<\min\{1,|abz/d|\}$. By the analytic
continuation, the restriction on $c$ can be relaxed. This completes
the proof of Theorem \ref{thm-a}.

%%%%%%%%%%%%%%%%%%%%%%%%%%%%%%%%%%%%%%%%%%%%%%%%%%%%%%%%%%%%%%%%%%%
\subsection{Proof of Theorem \ref{thm-b}}
%%%%%%%%%%%%%%%%%%%%%%%%%%%%%%%%%%%%%%%%%%%%%%%%%%%%%%%%%%%%%%%%%%%
%%%%%%%%%%%%%%%%%%%%%%%%%%%%%%%%%%%%%%%%%%%%%%%%%%%%%%%%%%%%%%%%%%%
By means of \eqref{bailey-a} and Watson's transformation formula for
three $_2\phi_1$ series (cf. \citu{gasper}{Appendix III.32}):
 \bnm
{_2\phi_1}\ffnk{cccccc}{q;z}{a,b}{c} &&\xqdn=\:
 \frac{(b,c/a,az,q/az;q)_{\infty}}{(c,b/a,z,q/z;q)_{\infty}}
 {_2\phi_1}\ffnk{cccccc}{q;\frac{qc}{abz}}{a,qa/c}{qa/b}
+idem(a;b),
 \enm
we attain
 \bnm
&&\xqdn{_2\psi_2}\ffnk{cccccc}{q;z}{a,b}{q^{1+m},d}\\
&&\xqdn\:=\:\frac{(a,b;q)_{-m}z^{-m}}{(q^{1+m},d;q)_{-m}}
\frac{(bq^{-m},d/a,azq^{-m},q^{1+m}/az;q)_{\infty}}{(dq^{-m},b/a,z,q/z;q)_{\infty}}\\
&&\:\:\times\:\:
{_2\phi_1}\ffnk{cccccc}{q;\frac{q^{1+m}d}{abz}}{aq^{-m},qa/d}{qa/b}+idem(a;b)\\
&&\xqdn\:=\:
\frac{(q,b,d/a,az,q/az,q^{1+m}/a;q)_{\infty}}{(q^{1+m},d,b/a,z,q/z,q/a;q)_{\infty}}
{_2\phi_1}\ffnk{cccccc}{q;\frac{q^{1+m}d}{abz}}{aq^{-m},qa/d}{qa/b}+idem(a;b).
 \enm
Split the ${_2\psi_2}$ series on the left hand side into two parts
to obtain
 \bmn
  \label{slater-a}
&&\xxqdn\sum_{k=0}^{\infty}\frac{(a,b;q)_k}{(q^{1+m},d;q)_k}z^k
+\sum_{k=1}^{\infty}\frac{(q/d;q)_k\prod_{i=1}^k(q^{1+m}-q^i)}{(q/a,q/b;q)_k}\bigg(\frac{d}{abz}\bigg)^k
 \nnm\\\nnm
&&\xxqdn\:=\frac{(q,b,d/a,az,q/az,q^{1+m}/a;q)_{\infty}}{(q^{1+m},d,b/a,z,q/z,q/a;q)_{\infty}}
\sum_{k=0}^{\infty}\frac{(qa/d;q)_k\prod_{i=1}^k(q^{1+m}-aq^i)}{(q,qa/b;q)_k}\bigg(\frac{d}{abz}\bigg)^k\\
&&\xxqdn\:+\:\:idem(a;b).
 \emn
 Define two functions $f(c)$ and $g(c)$ by
 \bnm
&&\xqdn f(c)=\sum_{k=0}^{\infty}\frac{(a,b;q)_k}{(c,d;q)_k}z^k
+\sum_{k=1}^{\infty}\frac{(q/d;q)_k\prod_{i=1}^k(c-q^i)}{(q/a,q/b;q)_k}\bigg(\frac{d}{abz}\bigg)^k,\\
&&\xqdn
g(c)=\frac{(q,b,d/a,az,q/az,c/a;q)_{\infty}}{(c,d,b/a,z,q/z,q/a;q)_{\infty}}
\sum_{k=0}^{\infty}\frac{(qa/d;q)_k\prod_{i=1}^k(c-aq^i)}{(q,qa/b;q)_k}\bigg(\frac{d}{abz}\bigg)^k\\
&&\:\:\:+\:\:idem(a;b).
 \enm
Then \eqref{slater-a} gives that
 \bmn \label{slater-b}
 f(c)=g(c)
\emn
 for $c=q^{1+m}$. In accordance with Lemma \ref{lemma}, \eqref{slater-b} is true for all
$|c|<1$. Through the analytic continuation, the restriction on $c$
could be relaxed. This finishes the proof of Theorem \ref{thm-b}.

%%%%%%%%%%%%%%%%%%%%%%%%%%%%%%%%%%%%%%%%%%%%%%%%%%%%%%%%%%%%%%%%%%%
\subsection{Proof of Theorem \ref{thm-c}}
%%%%%%%%%%%%%%%%%%%%%%%%%%%%%%%%%%%%%%%%%%%%%%%%%%%%%%%%%%%%%%%%%%%
%%%%%%%%%%%%%%%%%%%%%%%%%%%%%%%%%%%%%%%%%%%%%%%%%%%%%%%%%%%%%%%%%%%
In terms of \eqref{bailey-a} and the transformation formula
involving three $_2\phi_1$ series (cf. \citu{gasper}{Appendix
III.31}):
 \bnm
\quad{_2\phi_1}\ffnk{cccccc}{q;z}{a,b}{c} &&\xqdn=\:
 \frac{(abz/c,q/c;q)_{\infty}}{(az/c,q/a;q)_{\infty}}
 {_2\phi_1}\ffnk{cccccc}{q;\frac{qb}{c}}{c/a,qc/abz}{qc/az}
 \nnm\\&&\xqdn
 -\,\,\frac{(b,q/c,c/a,az/q,q^2/az;q)_{\infty}}{(c/q,qb/c,q/a,az/c,qc/az;q)_{\infty}}
 {_2\phi_1}\ffnk{cccccc}{q;z}{qa/c,\,qb/c}{q^2/c},
 \enm
we get
 \bnm
&&\xqdn{_2\psi_2}\ffnk{cccccc}{q;z}{a,b}{q^{1+m},d}\\
&&\xqdn\:=\:\frac{(a,b;q)_{-m}z^{-m}}{(q^{1+m},d;q)_{-m}}
\frac{(abzq^{-m}/d,q^{1+m}/d;q)_{\infty}}{(az/d,q^{1+m}/a;q)_{\infty}}
 {_2\phi_1}\ffnk{cccccc}{q;\frac{qb}{d}}{d/a,q^{1+m}d/abz}{qd/az}\\
&&\xqdn\:-\:\, \frac{(a,b;q)_{-m}z^{-m}}{(q^{1+m},d;q)_{-m}}
\frac{(bq^{-m},q^{1+m}/d,d/a,azq^{-m-1},q^{2+m}/az;q)_{\infty}}{(dq^{-m-1},qb/d,q^{1+m}/a,az/d,qd/az;q)_{\infty}}\\
&&\xqdn\:\times\:\:{_2\phi_1}\ffnk{cccccc}{q;z}{qa/d,\,qb/d}{q^{2+m}/d}\\
&&\xqdn\:=\:
\frac{(q,q^{1+m}/b,q/d,abz/d,qd/abz;q)_{\infty}}{(q/a,q/b,q^{1+m},az/d,q^{1+m}d/abz;q)_{\infty}}
{_2\phi_1}\ffnk{cccccc}{q;\frac{qb}{d}}{d/a,q^{1+m}d/abz}{qd/az}\\
 &&\xqdn\:-\:\,
\frac{(q,b,q/d,q^{2+m}/d,d/a,az/q,q^2/az;q)_{\infty}}{(q/a,q^{1+m},d/q,q^2/d,qb/d,az/d,qd/az;q)_{\infty}}
 {_2\phi_1}\ffnk{cccccc}{q;z}{qa/d,\,qb/d}{q^{2+m}/d}.
 \enm
Split the ${_2\psi_2}$-series on the left hand side into two parts
to gain
 \bmn
  \label{another-trans-a}
&&\xqdn\sum_{k=0}^{\infty}\frac{(a,b;q)_k}{(q^{1+m},d;q)_k}z^k
+\sum_{k=1}^{\infty}\frac{(q/d;q)_k\prod_{i=1}^k(q^{1+m}-q^i)}{(q/a,q/b;q)_k}\bigg(\frac{d}{abz}\bigg)^k
 \nnm\\\nnm
&&\xqdn\:=\:\frac{(q,q^{1+m}/b,q/d,abz/d,qd/abz;q)_{\infty}}{(q/a,q/b,q^{1+m},az/d,q^{1+m}d/abz;q)_{\infty}}
{_2\phi_1}\ffnk{cccccc}{q;\frac{qb}{d}}{d/a,q^{1+m}d/abz}{qd/az}\\
 &&\xqdn\:-\:\,
\frac{(q,b,q/d,q^{2+m}/d,d/a,az/q,q^2/az;q)_{\infty}}{(q/a,q^{1+m},d/q,q^2/d,qb/d,az/d,qd/az;q)_{\infty}}
 {_2\phi_1}\ffnk{cccccc}{q;z}{qa/d,\,qb/d}{q^{2+m}/d}.
 \emn
Define two functions $f(c)$ and $g(c)$ by
 \bnm
&&\qqdn\xqdn f(c)=\sum_{k=0}^{\infty}\frac{(a,b;q)_k}{(c,d;q)_k}z^k
+\sum_{k=1}^{\infty}\frac{(q/d;q)_k\prod_{i=1}^k(c-q^i)}{(q/a,q/b;q)_k}\bigg(\frac{d}{abz}\bigg)^k,\\
&&\qqdn\xqdn
g(c)=\frac{(q,c/b,q/d,abz/d,qd/abz;q)_{\infty}}{(q/a,q/b,c,az/d,cd/abz;q)_{\infty}}
{_2\phi_1}\ffnk{cccccc}{q;\frac{qb}{d}}{cd/abz,d/a}{qd/az}\\
 &&\qqdn\xqdn\qquad-\frac{(q,b,q/d,qc/d,d/a,az/q,q^2/az;q)_{\infty}}{(q/a,c,d/q,q^2/d,qb/d,az/d,qd/az;q)_{\infty}}
 {_2\phi_1}\ffnk{cccccc}{q;z}{qa/d,qb/d}{qc/d}.
 \enm
Then \eqref{another-trans-a} offers that
 \bmn \label{another-trans-b}
 f(c)=g(c)
\emn
 for $c=q^{1+m}$. According to Lemma \ref{lemma}, \eqref{another-trans-b} is right for all
$|c|<\min\{1,|d/q|,|abz/d|\}$. Via the analytic continuation, the
restriction on $c$ can by relaxed. This completes the proof of
Theorem \ref{thm-c}.

%%%%%%%%%%%%%%%%%%%%%%%%%%%%%%%%%%%%%%%%%%%%%%%%%%%%%%%%%%%%%%%%%%%
\section{Some related discuss}
%%%%%%%%%%%%%%%%%%%%%%%%%%%%%%%%%%%%%%%%%%%%%%%%%%%%%%%%%%%%%%%%%%%
The iteration of Theorem \ref{thm-a} produces another transformation
formula between two $_2\psi_2$ series due to Bailey \cito{bailey}:
\bnm
{_2\psi_2}\ffnk{cccccc}{q;z}{a,b}{c,d}=\frac{(az,bz,qc/abz,qd/abz;q)_{\infty}}{(q/a,q/b,c,d;q)_{\infty}}
{_2\psi_2}\ffnk{cccccc}{q;\frac{cd}{abz}}{abz/c,abz/d}{az,bz},
 \enm
where $\max\{|z|,|cd/abz|\}<1$.

Let $k$ denote the summation index of the $_2\psi_2$ series in
Theorem \ref{thm-b}. Replace $k$ by $-k$ to achieve
 \bnm
{_2\psi_2}\ffnk{cccccc}{q;\frac{cd}{abz}}{q/c,q/d}{q/a,q/b}
&&\xqdn\!=\frac{(q,b,c/a,d/a,az,q/az;q)_{\infty}}{(c,d,q/a,b/a,z,q/z;q)_{\infty}}
 {_2\phi_1}\ffnk{cccccc}{q;\frac{cd}{abz}}{qa/c,qa/d}{qa/b}\\
&&\xqdn\!+\:\:idem(a;b).
 \enm
Employing the substitutions $a\to q/c$, $b\to q/d$, $c\to q/a$,
$d\to q/b$, $z\to cd/abz$ in the last equation, we attain the case
$r=2$ of \citu{gasper}{Equation (5.4.5)}:
 \bnm
{_2\psi_2}\ffnk{cccccc}{q;z}{a,b}{c,d}
&&\xqdn\!=\frac{(q,q/d,c/a,c/b,abz/d,qd/abz;q)_{\infty}}{(q/a,q/b,c,c/d,cd/abz,qabz/cd;q)_{\infty}}
 {_2\phi_1}\ffnk{cccccc}{q;z}{qa/c,qb/c}{qd/c}\\
&&\xqdn\!+\:\:idem(c;d),
 \enm
provided $\max\{|z|, |cd/abz|\}<1$.

Performing the replacements $a\to b$, $b\to a$, $z\to cd/qab$ in
Theorem \ref{thm-c}, we have
 \bnm
{_2\psi_2}\ffnk{cccccc}{q;\frac{cd}{qab}}{a,b}{c,d} &&\xqdn=\:
\frac{(c/a,c/q,q^2/c,q/d;q)_{\infty}}{(c,c/qa,q/a,q/b;q)_{\infty}}
{_2\phi_1}\ffnk{cccccc}{q;\frac{qa}{d}}{d/b,q}{q^2a/c}
 \\&&\xqdn
 -\,\,\frac{(q,a,d/b,q/d,qc/d,cd/q^2a,q^3a/cd;q)_{\infty}}{(q/b,c,d/q,q^2/d,qa/d,c/qa,q^2a/c;q)_{\infty}}
\\&&\xqdn
\times\,\,{_2\phi_1}\ffnk{cccccc}{q;\frac{cd}{qab}}{qa/d,qb/d}{qc/d}.
 \enm
Calculating the second ${_2\phi_1}$ series on the right hand side by
$q$-Gauss summation formula (cf. \citu{gasper}{Appendix II.8}):
 \bnm
 {_2\phi_1}\ffnk{cccccc}{q;\frac{c}{ab}}{a,b}{c}=\frac{(c/a,c/b;q)_{\infty}}{(c,c/ab;q)_{\infty}},
 \enm
we obtain
 \bnm
{_2\psi_2}\ffnk{cccccc}{q;\frac{cd}{qab}}{a,b}{c,d} &&\xqdn=\:
\frac{(c/a,c/q,q^2/c,q/d;q)_{\infty}}{(c,c/qa,q/a,q/b;q)_{\infty}}
{_2\phi_1}\ffnk{cccccc}{q;\frac{qa}{d}}{d/b,q}{q^2a/c}
 \\&&\xqdn
 +\,\,\ffnk{cccccccc}{q}{q,a,c/b,d/b,cd/qa,q^2a/cd}{c,d,q/b,qa/c,qa/d,cd/qab}_{\infty}.
 \enm
Applying another form of Heine's transformation (cf.
\citu{gasper}{Appendix III.1})
\[{_2\phi_1}\ffnk{cccccc}{q;z}{a,b}{c}=\frac{(b,az;q)_{\infty}}{(c,z;q)_{\infty}}{_2\phi_1}\ffnk{cccccc}{q;b}{c/b,z}{az}\]
to the last equation, we recover the following formula due to Chu
\cito{chu-b}:
 \bnm
{_2\psi_2}\ffnk{cccccc}{q;\frac{cd}{qab}}{a,b}{c,d}
&&\xqdn=\:\ffnk{cccccccc}{q}{q,a,c/b,d/b,cd/qa,q^2a/cd}{c,d,q/b,qa/c,qa/d,cd/qab}_{\infty}
 \\&&\xqdn
 +\:\:a\ffnk{cccccccc}{q}{q,qa/b,q/c,q/d}{qa/c,qa/d,q/a,q/b}_{\infty}\qdn
 {_2\phi_1}\ffnk{cccccc}{q;q}{qa/c,qa/d}{qa/b},
 \enm
where the convergent condition is $|cd/qab|<1$.

\textbf{Acknowledgments}

 The work is supported by the National Natural Science Foundation of China (No. 11661032).

%%%%%%%%%%%%%%%%%%%%%%%%%%%%%%%%%%%%%%%%%%%%%%%%%%%%%%%%%%%%%%%%%%%

%%%%%%%%%%%%%%%%%%%%%%%%%%%%%%%%%%%%%%%%%%%%%%%%%%%%%%%%%%%%%%%%%%%
%%%%%%%%%%%%%%%%%%%%%%%%%%%%%%%%%%%%%%%%%%%%%%%%%%%%%%%%%%%%%%%%%%%
%%%%%%%%%%%%%%%%%%%%%%%%%%%%%%%%%%%%%%%%%%%%%%%%%%%%%%%%%%%%%%%%%%%

\end{document}